\documentclass[reqno]{amsart}

\usepackage{amsmath}
\usepackage{amssymb}
\usepackage{graphicx}
\usepackage{color}
\usepackage[all]{xy}
\usepackage{tikz-cd}
\usepackage{amsthm}

\theoremstyle{plain}

\newtheorem{para}{}[section]
\newtheorem{thm}[para]{Theorem}
\newtheorem{prop}[para]{Proposition}
\newtheorem{lemma}[para]{Lemma}
\newtheorem{cor}[para]{Corollary}

\newtheorem*{quest}{Question}

\theoremstyle{remark}

\newtheorem{remark}[para]{Remark}
\newtheorem{remarks}{Remarks}

\theoremstyle{definition}

\newcommand{\wh}[1]{\widehat{#1}}
\newcommand{\co}{\colon\thinspace}

\newcommand{\bound}{\partial}
\newcommand{\pie}{\pi_1}
\newcommand{\normal}[1]{\langle \! \langle #1 \rangle \! \rangle}

\newcommand{\C}{\mathbb{C}}
\renewcommand{\H}{\mathbb{H}}
\newcommand{\Q}{\mathbb{Q}}
\renewcommand{\O}{\mathbb{O}}
\newcommand{\R}{\mathbb{R}}
\renewcommand{\S}{\mathbb{S}}
\newcommand{\T}{\mathbb{T}}
\newcommand{\Z}{\mathbb{Z}}

\newcommand{\calF}{\mathcal{F}}
\newcommand{\calG}{\mathcal{G}}

\newcommand{\calT}{\mathcal{T}}
\newcommand{\calU}{\mathcal{U}}

\newcommand{\sfC}{{\sf C}}
\newcommand{\sfD}{{\sf D}}
\newcommand{\sfX}{{\sf X}}

\newcommand{\tr}{\mathrm{Tr}}

\newcommand{\PSL}{A point $\chi=(X_0,Z_0) \in \C^2$ is an irreducible character in $\overline{\sfD}_\alpha$  if and only if $X_0Z_0 \neq 0$ and $\chi$ satisfies the polynomial $\calT_0(\alpha)$.}

\newcommand{\Factors}{If $\alpha' \in \wh{G}_\alpha \{ \alpha, \infty\}$ then every factor of $\calT(\alpha)$ divides $\calT(\alpha')$.}

\newcommand{\Epimorphism}{Suppose that $\alpha \in \wh{\Q}$ and $M(\alpha)$ is hyperbolic.  If $\alpha' \in \wh{G}_\alpha\{\alpha, \infty\}$ then there is an epimorphism $\Gamma_\alpha \to \Gamma_{\alpha'}$ taking $k_j\in \Gamma_\alpha$ to $k_j \in \Gamma_{\alpha'}$.}

\newcommand{\Riley}{ If $\alpha \in \Q_0$ has odd denominator then, up to sign,  $\Lambda_\alpha$ is the Riley polynomial.}

\newcommand{\TO}{For $p/q \in \wh{\Q}$, $\calT_0(p/q)  \in \Z[X, Z]$.}

\begin{document}

\title{Farey recursion and the character varieties for 2-bridge knots} 

\author{Eric Chesebro}
\address{Department of Mathematical Sciences, University of Montana} 
\email{Eric.Chesebro@mso.umt.edu} 

\begin{abstract}  
We describe the $\text{(P)SL}_2 \C$ character varieties of all 2-bridge knots and the diagonal character varieties for all 2-bridge links in terms of a set of polynomials defined using {\em Farey recursion}.  \end{abstract}

\maketitle

\section{Introduction}

When $M$ is a hyperbolic 3-manifold with a single cusp, the $\text{(P)SL}_2 \C$ character variety for $M$ has proven to be a valuable tool in study the topology and geometry of $M$.  There are too many examples to list here, but some classics are  \cite{CS}, \cite{CGLS}, and \cite{Th_notes}.  Although the theory of character varieties is well-developed and very useful, it is notoriously difficult to compute character varieties in specific examples, particularly as the complexity of the manifolds increases.

A large class of link complements in which some computational progress has been made is the class of 2-bridge knots.  Independent of the study of character varieties, these attractive links have a long history of mathematical interest.  A highlight in this history, is the connection between 2-bridge links, continued fractions, and the Farey graph.  In 1956, Schubert showed that 2-bridge links are classified by rational numbers and showed that continued fraction expansions of rational numbers give a natural way to diagram their corresponding link as the closure of a 4-strand braid.  Famously in \cite{HT}, Hatcher and Thurston used this classification and its relationship to the Farey graph to classify the essential surfaces in 2-bridge link complements.  Sakuma and Weeks \cite{SW} also used the relationship between these links and the Farey graph to construct beautiful ideal triangulations of the complements of 2-bridge links.  Independently, \cite{PTGTKG} and \cite{GF} showed later that these triangulations are geometric and canonical in the sense of \cite{EP}.

In \cite{RR} and \cite{RR2}, Riley began a thorough study of the $\text{(P)SL}_2 \C$ representations of 2-bridge knot groups.  Among other things, he introduced the Riley polynomial for a 2-bridge knots, whose solutions correspond representations which take peripheral elements to parabolic matrices.  Riley's formulas have been used in \cite{HLM}, \cite{HS3}, \cite{MPV}, and \cite{PT} to collectively give recursive formulas for the $\text{(P)SL}_2 \C$ character varieties of 2-bridge links obtained by doing $1/n$ Dehn filling on unknotted components of the links $6_1^3$, $6_2^3$, $6_3^2$, $8_2^4$, and $8_9^3$.  This is a large class of links with arbitrarily large crossing numbers.  However, being Dehn fillings on a hand full of links, it is quick to check that the volumes for these examples are bounded by 7.4, while the volumes of 2-bridge links can be arbitrarily large.  Nonetheless, including Riley's original results, these formulas have enjoyed a variety of applications including \cite{Bu},  \cite{Chu}, \cite{G}, \cite{HS3}, \cite{HS2},  \cite{HS1}, and \cite{RW}.

Here, we study the $\text{(P)SL}_2 \C$ character varieties of all 2-bridge knots together with the {\em diagonal} character varieties of the 2-bridge links.  This makes it possible to see an attractive relationship between the defining polynomials for all of these algebraic sets.  We refer to this relationship as {\em Farey recursion}.

By identifying rational numbers $p/q$ with the point $p/q+i/q$ in the upper half $\H^2$ of $\C$, we can connect the rationals in a collection of interlocking triangles to form the Stern-Brocot diagram, see Figures \ref{fig: calG} and \ref{fig: Tri}.  Each of these triangles contains exactly one rational number $p/q+i/q$ in its interior, which we call its {\em center}.  The vertices that occur along the edges of this triangle form a bi-infinite sequence of rational numbers, which we refer to as the {\em boundary sequence} for $p/q$.  In Section \ref{sec: Farey recursion}, we argue that there is a unique function \[\calT \co \Q\cup \{\infty\} \to \Z[x,z]\] 
with the following properties.
\begin{enumerate}
\item $\calT(0)=x$, $\calT(\infty)=0$, and $\calT(1)=z$;
\item Given $p/q \in \Q$, and three successive terms $\gamma_1, \gamma_2, \gamma_3$ in the boundary sequence for $p/q$, 
\[ \calT(\gamma_3) = -\calT(\gamma_1)+\calT(p/q) \, \calT(\gamma_2).\]
\end{enumerate}
Property (2) is the {\em Farey recursion} condition.  We can imagine $\calT(\Q)$ as a set of polynomials that are ``linearly'' recursive on the Stern-Brocot diagram, rather than on a line.  Property (2) also makes it possible to compute the values of $\calT$ relatively efficiently with a computer.  Lastly, property (2) often allows us to make elementary arguments regarding $\calT$.  For instance, if we define $f \co \Q \cup \{ \infty\} \to \Z[x,z]$ by 
\[ f(p/q) = \begin{cases} x & \text{if } p \text{ is even and } q \text{ is odd}\\
z & \text{if } p \text{ and } q \text{ are odd}\\
xz & \text{if } p \text{ is odd and } q \text{ is even}\end{cases}\]
and set $X=x^2$, and $Z=z^2$, an elementary inductive argument  (Lemma \ref{lem: fac}) shows that 
\[ \calT_0(\alpha)= \frac{\calT(\alpha)}{f(\alpha)} \in \Z[X,Z]\]
for every $\alpha \in \Q \cup \{ \infty\}$.

The main theorem of this paper shows that the $\text{PSL}_2 \C$ character variety for the 2-bridge link $L(\alpha)$ associated to $\alpha$ is naturally identified with the complex affine algebraic set determined by $\calT_0(\alpha)$.

 \theoremstyle{plain}
\newtheorem*{PSL thm}{Theorem \ref{thm: PSL}}
\begin{PSL thm} \PSL \end{PSL thm}

Here $\overline{\sfD}_\alpha$ is the {\em diagonal} $\text{PSL}_2 \C$ character variety for the link $L(\alpha)$.    In the case that $L(\alpha)$ is a knot, $\overline{\sfD}_\alpha$ is the full character variety of irreducible characters.

It is relatively efficient to compute the character variety polynomials $\calT_0(\alpha)$.  With a simple program implemented in CoCalc \cite{sage}, we are able to compute all 6079 values for $\calT_0$ on the set
\[ \left\{ \frac{p}{q} \in \Q \cap [0,1/2] \, \bigg| \, q<200 \right\}\]
in roughly 8 minutes. 

The main theorem of \cite{ORS} gives conditions on $\alpha'$ which guarantee existence of an epimorphism from the link group of $L(\alpha)$ to that of $L(\alpha')$.  In Subsection \ref{subsec: epi}, we employ a natural action of the modular group on the set of polynomials $\calT(\Q \cup \{ \infty\})$ to understand when the factors of $\calT_0(\alpha)$ must divide $\calT_0(\alpha')$, see  Corollary \ref{cor: factors}.  Ultimately, this provides an elementary argument which reproduces the result in \cite{ORS}, see  Corollary \ref{cor: epi}.  

In Subsection \ref{subsec: MPV}, we describe the relationship between our formulas and the formulas in \cite{MPV}.
 
In \cite{RR}, Riley defines a polynomial in $\Z[X]$ for each rational number with odd denominator.  These polynomials have since been refered to as {\em Riley polynomials}.  Among other things, Riley proved that his polynomials are monic and their roots give rise to normalized {\em p-reps} for the corresponding 2-bridge knot groups.  We discuss this in Section \ref{sec: Riley}, where we also show that by setting $Z=-X$ in $\calT_0(\alpha)$ we recover the Riley polynomial and thus extend the definition for Riley polynomials to rational numbers with even denominators.  Our approach provides an efficient way to compute Riley polynomials as well as access to straightforward inductive arguments concerning their properties.  

Again using CoCalc \cite{sage}, we are able to compute all $\Lambda_\alpha$ for each the the 13,662 numbers in 
\[ \left\{ \frac{p}{q} \in \Q \cap [0,1/2] \, \bigg| \, q<300 \right\}\]
in roughly 12 minutes. 

To demonstrate the inductive facility of our setup, we give a quick proof that $\Lambda_\alpha$ is monic for every $\alpha$.  

\subsection*{Acknowledgement} Thanks to Kelly McKinnie for her numerous corrections and suggestions.

\section{Surfaces} The material in this section is standard.   Roughly, we follow the approach in \cite{Ford}.  Let $\mathbb{X}$ be the punctured plane $\C-\Z[i]$ equipped with its usual Euclidean geometry.  Define subgroups $\Gamma_O,\Gamma_S, \Gamma_T \subset\text{Isom}(\mathbb{X})$ by  \begin{align*}
\Gamma_O &= \left\langle \text{ order-2 rotations centered at the points in } \frac{1}{2}\,  \Z[i] \right\rangle \\
\Gamma_S &= \Big\langle \text{ order-2 rotations centered at the points in } \Z[i] \, \Big\rangle \\
\Gamma_T &= \big\langle \text{  horizontal and vertical unit translations } \big\rangle.
\end{align*}  Then $\O=\mathbb{X}/\Gamma_O$ is a $(2,2,2,\infty)$-pillowcase orbifold, $\S=\mathbb{X}/\Gamma_S$ is a 4-punctured sphere, and $\T=\mathbb{X}/\Gamma_T$ is a once punctured torus.  We also have a commutative diagram of covers.  \[ \xymatrix{  & \mathbb{X} \ar[dl] \ar[dr] \\ \T  \ar[dr]_2 &  & \S \ar[dl]^4 \\ & \O} \]

Take $\tilde{x}= (1/4, 1/4)$ as a basepoint in $\mathbb{X}$.  Consider the paths  \begin{align*}
&\frac{e^{\frac{3 \pi i }{4}}}{2\sqrt{2}} \, e^{t \pi i  } + \frac{1}{2} &
&\frac{e^{\frac{5 \pi i }{4}}}{2\sqrt{2}} \, e^{t \pi i  } + \frac{1}{2}(1+i) &
&\frac{e^{\frac{7 \pi i }{4}}}{2\sqrt{2}} \, e^{t \pi i  } + \frac{i}{2} 
\end{align*} in $\mathbb{X}$ based at $\tilde{x}$ and let $\tilde{p}$, $\tilde{q}$, and $\tilde{r}$ be their respective homotopy classes (rel $\bound$).  Then, if $p$, $q$, and $r$ are the corresponding elements of $\pie(\O)$,   \[ \pie(\O) = \left\langle p, q, r  \mid  p^2, q^2, r^2  \right\rangle.\]   Let $k_0=rpq$ and observe that $k_0$ is primitive and peripheral.

Define \begin{align*} a&=rq & &\text{and} & b&=qp. \end{align*} We also let $A=a^{-1}$ and $B=b^{-1}$.  As a subgroup of $\pie(\O)$, $\pie(\T)$ is the rank-2 free group generated by $a$ and $b$.  Moreover, the corresponding deck transformations are translations by $-1$ and $i$ respectively.  The commutator $[a,b]=k_0^2$ is primitive and peripheral in $\pie \T$.  

\begin{figure}[h] 
   \centering
   \includegraphics[width=2in]{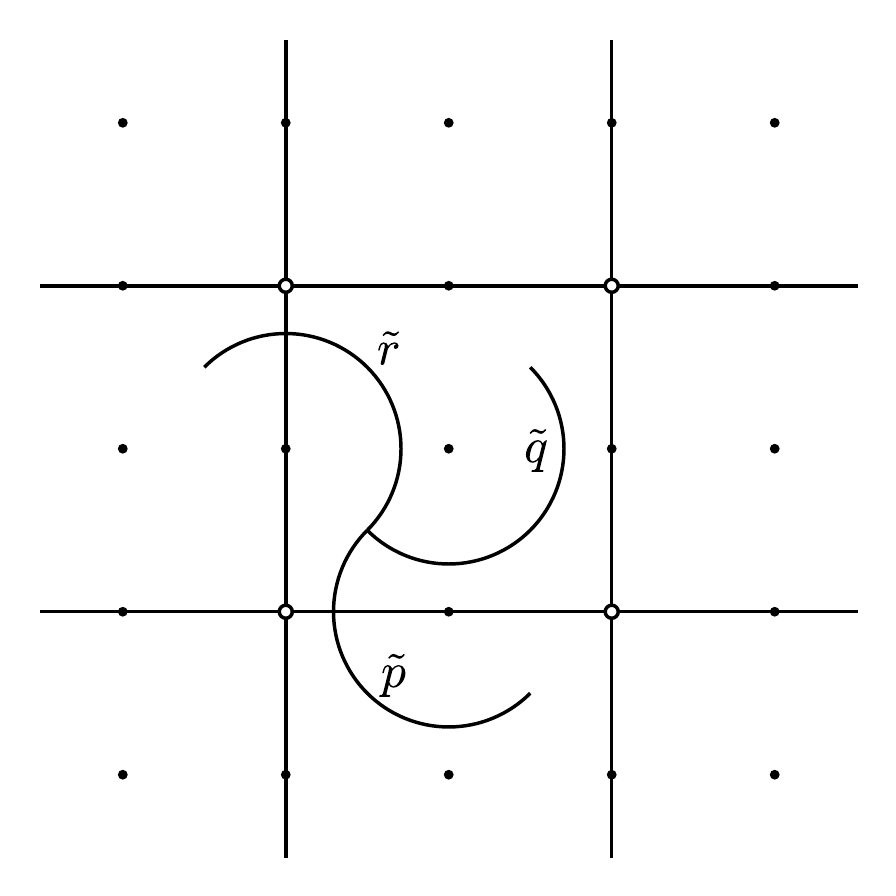} 
   \caption{Paths in $\mathbb{X}$.}   \label{fig: X}
\end{figure}

Define distinguished conjugates \begin{align*} k_1&= pk_0p & k_2&=qk_0q & k_3&=rk_0r \end{align*} of $k_0$.  As a subgroup of $\pie(\O)$, $\pie(\S)$ is the rank-3 free group generated by the peripheral elements $k_0$, $k_1$, and $k_2$. The product $k_0k_2k_1k_3$ is trivial and \begin{align*}  a^2&= k_3^{-1}k_2^{-1} & B^2&= k_3k_0 & (aB)^2&=k_3^{-1}k_2k_3k_0. \end{align*}

\section{$\mathbf{ \text{PSL}_2 \C}$ representations} \label{sec: reps} Given ${\bf z}=(x,z) \in \C^2$ we can find $t \in \C-\{0\}$ so that \begin{align} \label{eq: Markov} x^2+z^2+\left( t - t^{-1} \right)^2 &=0. \end{align} We refer to the pair ${\bf z},t$ as a {\bf parameter pair}.  Every parameter pair ${\bf z},t$ determines a representation $\rho_{{\bf z}, t} \co \pie \O \to \text{PSL}_2 \C$ as follows. If $x=0$, define
\begin{align*}
\rho_{{\bf z}, t}(p) & = \begin{bmatrix}  0&1\\-1&0 \end{bmatrix} &
\rho_{{\bf z}, t}(q) & =  \begin{bmatrix} i&0 \\ 0 & -i \end{bmatrix} &
\rho_{{\bf z}, t}(r) & =  \begin{bmatrix} 0& it \\ \frac{i}{t} & 0 \end{bmatrix}.
\end{align*}
If $x \neq 0$ and $z=0$, define
\begin{align*}
\rho_{{\bf z}, t}(p) & = \begin{bmatrix} i&0 \\ 0 & -i \end{bmatrix} &
\rho_{{\bf z}, t}(q) & =  \begin{bmatrix}  0&1\\-1&0 \end{bmatrix} &
\rho_{{\bf z}, t}(r) & =  \begin{bmatrix} 0& it \\ \frac{i}{t} & 0 \end{bmatrix}.
\end{align*}
If $xz \neq 0$, define
\begin{align*}
\rho_{{\bf z}, t}(p) & =  \frac{1}{x} \, \begin{bmatrix} t^{-1}-t & -1 \\ -z^2 & t-t^{-1} \end{bmatrix} \\
\rho_{{\bf z}, t}(q) & = \frac{1}{z} \,  \begin{bmatrix} 0 & 1 \\ -z^2 & 0 \end{bmatrix} \\
\rho_{{\bf z}, t}(r) & =  \frac{1}{x} \, \begin{bmatrix} -z \,t^{-1} & z+t^2 \,z^{-1} -z^{-1} \\ z\, (1-t^{-2}) & z \, t^{-1} \end{bmatrix}.
\end{align*}

Straightforward calculations verify that, in all cases,
\begin{align*}
\tr \big( \rho_{{\bf z}, t}(a) \big) &= \pm x & \tr \big( \rho_{{\bf z}, t}(b) \big) &= 0 & \tr \big( \rho_{{\bf z}, t}(aB) \big) &= \pm z.
\end{align*}
Also, if $x=0$, then 
\[ \rho_{{\bf z}, t}(k_0) = \rho_{{\bf z}, t}(k_1^{-1})= \begin{bmatrix} t & 0 \\ 0 & t^{-1} \end{bmatrix}, \]
if $x \neq 0$ and $z=0$, then 
\[ \rho_{{\bf z}, t}(k_0) = \rho_{{\bf z}, t}(k_1)= \begin{bmatrix} t & 0 \\ 0 & t^{-1} \end{bmatrix}, \]
and, if $x z \neq 0$, then
\begin{align*}
\rho_{{\bf z}, t}(k_0) &= \begin{bmatrix} t & 1 \\ 0 & t^{-1} \end{bmatrix} & \rho_{{\bf z},t}(k_1) &= \begin{bmatrix} t & 0 \\ z^2 & t^{-1} \end{bmatrix}. 
\end{align*}
Note that $\rho_{{\bf z}, t}$ is irreducible if and only if $xz \neq 0$.

\section{Slopes and essential simple closed curves}
Lines in $\mathbb{X}$ with rational slopes descend to essential simple closed curves on $\O$, $\T$, and $\S$.  In fact, there are bijections from $\wh{\Q}=\Q \cup \{ 1/0 \}$ to the set of conjugacy classes in $\pie(\O)$, $\pie(\T)$, and $\pie(\S)$ whose corresponding free homotopy classes are represented by essential simple closed curves.  For example, in $\pie(\T)$ the slope of the class containing $a$ is zero and the slope of the class containing $b$ is infinite.  

Given $\alpha \in \Q_0= \Q \cap [0,1]$ and an $\text{SL}_2 \C$ representation of $\pie(\O)$, $\pie(\S)$, or $\pie(\T)$, it is surprisingly easy to find elements of $\pie(\S)$ and $\pie(\T)$ which represent the slope $\alpha$ and to compute their traces under the representation.  To show how this is done, it will be helpful to first discuss the Farey sum and the Stern-Brocot diagram.

\subsection{Farey sums and the Stern-Brocot diagram} \label{ssec: FSSBD}
Throughout the remainder of the paper, when we write elements of $\Q$ as integer quotients, we will always write them in lowest terms with positive denominators.  If $p/q$ and $r/s$ are elements of $\wh{\Q}$ such that $ps-rq=\pm1$, then $\{ p/q, r/s \}$ is called a {\bf Farey pair}.   The {\bf Farey sum} of $p/q$ and $r/s$ is defined as \[ \frac{p}{q} \oplus \frac{r}{s} = \frac{p+r}{q+s}.\]  It is remarkable that if $\{ p/q, r/s \}$ is a Farey pair then $p+r$ and $q+s$ are  relatively prime.  Notice that every Farey pair is a subset of exactly two {\bf Farey triples}.

\begin{figure}[h] 
   \centering
   \includegraphics[width=5in]{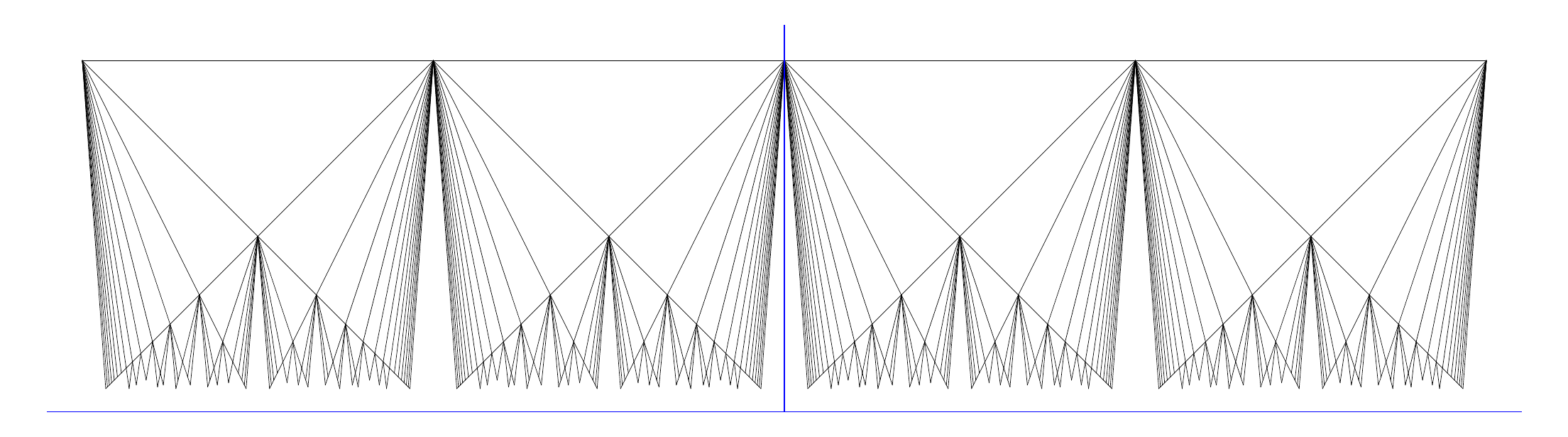} 
   \caption{The portion of $\calG$ spanned by fractions in $[-2,2]$ with denominator at most 15.}   \label{fig: calG}
\end{figure}

There is a closely related geometric graph $\calG$ embedded in the upper-half space $\H^2 \subset \C$.  The vertices of $\calG$ correspond to $\Q$ and the edges correspond to Farey pairs.  More precisely, for $\frac{p}{q} \in \Q$, define $v_{\frac{p}{q}} = \frac{p}{q}+\frac{i}{q} \in \H^2$.  The set \[ \left\{ v_\alpha \, | \, \alpha \in \Q \right\}\] constitutes the vertices of $\calG$ and the edges of $\calG$ are the straight line segments between Farey pairs of vertices.  (We will often blur the distinction between a number $\alpha \in \Q$ and its vertex $v_\alpha$.)  See Figure \ref{fig: calG} below as well as Figure 1.3 of \cite{HO}.  Here, we refer to $\calG$ as the {\bf Stern-Brocot diagram}.  Elementary properties of the Farey sum and Farey pairs imply that every intersection between a distinct pair of edges in $\calG$ occurs at a vertex of $\calG$.  

\begin{figure}[h] 
   \centering
   \includegraphics[width=1.2in]{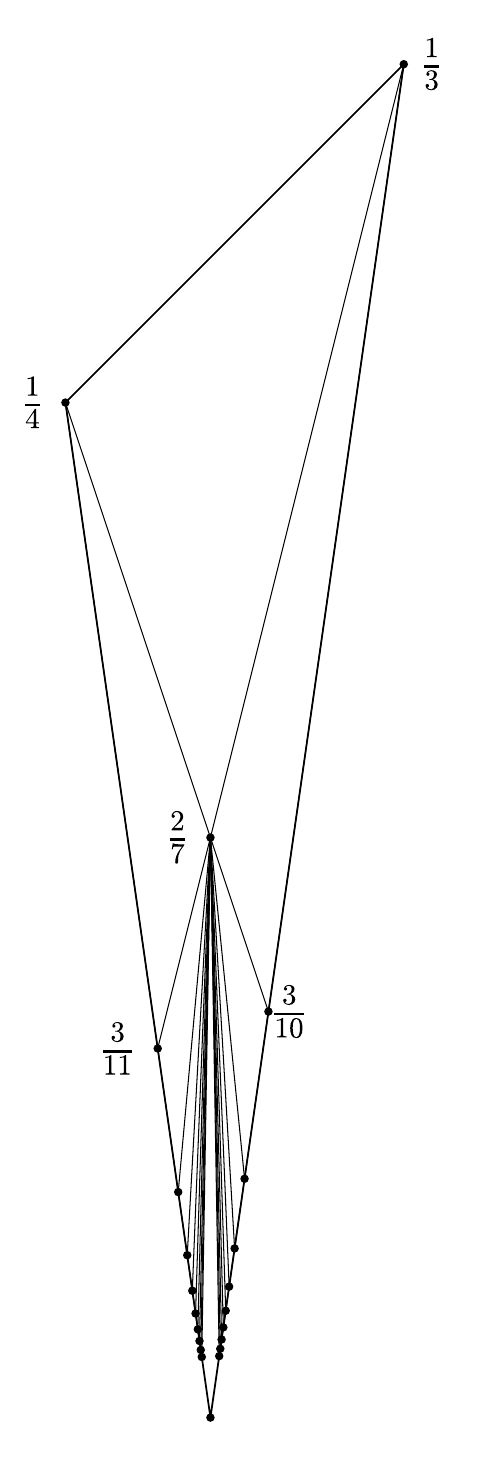} 
   \caption{$\Delta\left( \frac{2}{7} \right)$.}   \label{fig: Tri}
\end{figure}

Suppose $\alpha \in \Q$.  The graph theoretical neighborhood $\Delta(\alpha)$ of $v_\alpha$ in $\calG$ is, by definition, the subgraph spanned by the vertices which are Farey pairs with $\alpha$, see Figure \ref{fig: Tri}.

If $\alpha \in \Q - \Z$ then $\Delta(\alpha)$ corresponds to a triangle with an ideal vertex at $\alpha \in \R$.  This triangle contains exactly one vertex of $\calG$ in its interior, namely $v_\alpha$.  For this reason, we refer to $\Delta(\alpha)$ as the {\bf triangle centered at} $\mathbf{v_\alpha}$.  Its ideal vertex $\alpha$ is called its {\bf tip}.  Emanating from the tip of $\Delta(\alpha)$ are its left and right sides whose slopes are $\pm q$, where $q$ is the denominator of $\alpha$.  The vertices at the tops of these edges are called the left and right {\bf corners} of $\Delta(\alpha)$.  The {\bf positive} (resp. negative) {\bf vertex sequence} $\bound_+(\alpha)$ (resp. $\bound_-(\alpha)$) of $\Delta(\alpha)$ is defined to be the ordered sequence of vertices down the side of $\Delta(\alpha)$ whose slope is positive (resp. negative).   In fact, if $\kappa_\pm$ is the corner on this side, then
\begin{align} \label{bound} \bound_\pm ( \alpha) = \left\{\left( \oplus_0^j \alpha \right) \oplus \kappa_\pm   \right\}_{j=0}^\infty. \end{align}
We write $\bound(\alpha)$ to denote the bi-infinite sequence formed by reversing $\bound_-(\alpha)$ and juxaposing the sequence $\bound_+(\alpha)$ and we refer to $\bound (\alpha)$ as the {\bf boundary sequence} for $\alpha$.

The Stern-Brocot diagram and our notions of triangles, edges, and corners all extend nicely to include an additional vertex for $1/0$.   To do this, add a vertex $v_{1/0}$ and an edge from $v_{1/0}$ to each integer vertex of $\calG$.  The triangles centered at integers and at $1/0$ involve the new vertex $v_{1/0}$.  They are exceptional in that they each have two sides but only one corner.  We explicitly define their vertex sequences by listing their terms.  Let $n \in \Z$ and define
\begin{align*}
\bound_- \left(1/0\right) &= \left\{ 0, -1, -2, -3, \ldots \right\} \\
\bound_+ \left(1/0\right) &= \left\{ 0, 1, 2, 3, \ldots \right\} \\
\bound_- (n) &= \left\{ \frac{1}{0}, n-1, \frac{2n-1}{2}, \frac{3n-1}{3},  \ldots \right\} \\
\bound_+ (n) &= \left\{ \frac{1}{0}, n+1, \frac{2n+1}{2}, \frac{3n+1}{3},  \ldots \right\}.
\end{align*}
In the boundary sequences $\bound(1/0)$ and $\bound(n)$, we do not repeat the duplicate corner, for instance, $\bound(1/0)$ is the ordered sequence of integers.  Also, if we allow $\frac{1}{0}=\frac{-1}{0}$ and $\frac{0}{1}=\frac{0}{-1}$ and choose appropriate representatives, then the Farey sum formulas (\ref{bound}) still apply in these exceptional cases.  

There is a direct correlation between $\calG$ and the Farey Graph.  In \cite{Hat}, Hatcher points out that if the vertices of $\calG$ are pushed straight down to $\R \subset \H^2$ and the edges follow their endpoints, ultimately forming geodesics in $\H^2$, the result is the Farey graph.

The next lemma will be useful later in this paper.

\begin{lemma} \label{lem: properside}
Given $\alpha \in \Q-\Z$ there is a number $\gamma \in \Q$ where $\alpha$ is at least the third term of a vertex sequence for $\Delta(\gamma)$.
\end{lemma}

\proof
If the denominator of $\alpha$ is two, then $\alpha$ is the third term of a sequence $\bound_+(n)$ for some $n \in \Z$.  So, we may assume that the denominator of $\alpha$ is larger than two.

Since the denominator of $\alpha$ is larger than two, $\Delta(\alpha)$ has two corners, $\kappa_\pm \in \Q$, each of which lies above $\alpha$ in $\calG$.  The Farey pair condition forces the denominators of $\kappa_+$ and $\kappa_-$ to be distinct.  Without loss of generality,  assume that the denominator of $\kappa_-$ is smaller than that of $\kappa_+$.  

We know that $\alpha$ comes after $\kappa_+$ in the sequence $\bound_+(\kappa_-)$.  Because the denominator of $\kappa_+$ is larger than that of $\kappa_-$, we also know that $\kappa_+$ is not the first term of this sequence.  Hence, $\alpha$ is at least the third term of $\bound_+(\kappa_-)$.
\endproof

Let $G$ be the extended modular group 
\[ G=\text{PSL}_2 \Z \oplus \left\langle \begin{bmatrix} -1&0 \\ 0&1 \end{bmatrix} \right\rangle.\]
Every element of $G$ determines an invertible M\"{o}bius transformation and these transformations preserve Farey pairs when applied to $\wh{\Q}$.  This gives an action of $G$ on $\calG$.  Although this action preserves triangles, it does not preserve their corners.  For instance, if $\phi=\left[ \begin{smallmatrix} -1&1 \\ -4&3 \end{smallmatrix} \right]$ then $\phi$ fixes $\frac{1}{2}$ and preserves $\bound\left( \frac{1}{2} \right)$.  On the other hand, $\phi$ acts non-trivially on $\bound \left( \frac{1}{2} \right)$ taking each term to its neighbor on the left.

\subsection{Words of slope $\mathbf{\alpha}$}  There is a function $\Omega \co \wh{\Q}_0 \to \pie(\O)$ which is uniquely determined by the following two properties
\begin{enumerate}
\item If $\{\alpha, \beta\}$ is a Farey pair in $\wh{\Q}_0$ and $\alpha<\beta$ then
\[ \Omega(\alpha \oplus \beta) = \Omega(\alpha) \, \Omega(\beta).\]
\item $ \Omega(0)= a $ and $ \Omega(1/0)=B$.
\end{enumerate}
Here, we use the convention that $0/1<1/0$.  

If $\alpha$ is a positive rational number which is not in $\Q_0$, we define $\Omega(\alpha)$ to be the word obtained from $\Omega(1/\alpha)$ by swapping $a$'s and $B$'s.  If $\alpha$ is a negative rational number, we define $\Omega(\alpha)$ to be the word obtained from $\Omega(-\alpha)$ by replacing $B$'s with $b$'s.  This determines an extension $\Omega \co \wh{\Q} \to \pie(\O)$.  The next proposition is established in \cite{MDS} and \cite{S}.  

\begin{prop}
$\Omega(\alpha)$ has slope $\alpha$ in $\O$ and $\T$. 
\end{prop}

\begin{remark}If we consider the cover $\S \to \O$, we see that $\Omega(\alpha)^2$ has slope $\alpha$ in $\S$. \end{remark}

\section{Trace functions and Farey recursion.} \label{sec: Farey recursion}

\subsection{Farey recursion}  
Suppose that $R$ is a commutative ring and $d \co \wh{\Q} \to R$ is a function.  A function $\calF \co \wh{\Q} \to R$ is a {\bf Farey recursive function} (FRF) with {\bf determinant} $d$ if, whenever $\alpha$ and $\gamma$ make a Farey pair, then 
\begin{align} \label{eq: FRF} \calF(\alpha \oplus \alpha \oplus \gamma) &=  - d(\alpha) \, \calF(\gamma) +  \calF(\alpha) \, \calF(\alpha \oplus \gamma).\end{align}

Notice that if $\{ \alpha, \gamma \}$ is a Farey pair then $\{ \alpha, \gamma, \alpha \oplus \gamma, \alpha \oplus \alpha \oplus \gamma \}$ are the vertices of a pair of adjacent triangles in the Farey graph.  

Suppose that $\calF$ is an FRF with determinant $d$ and $\alpha \in \wh{\Q}$.  Choose $\epsilon \in \{ \pm 1\}$ and let $a_j \in R$ be the element obtained by applying $\calF$ to the $j^\text{th}$ element of $\bound_\epsilon(\alpha)$.
Then the sequence $\{ a_j \}_0^\infty$ is linearly recursive and satisfies
\begin{align} \label{eq: mtx} \begin{pmatrix} 0 & 1 \\ -d(\alpha) & \calF(\alpha) \end{pmatrix} \ \begin{pmatrix} a_{j-1} \\ a_j \end{pmatrix} = \begin{pmatrix} a_j \\ a_{j+1} \end{pmatrix}\end{align} for every $j\geq 1$.  We call the $2 \times 2$ matrix in Equation \ref{eq: mtx} the {\bf recursion matrix} at $\alpha$ for $\calF$.

Henceforth, we will assume that $R$ is a ring with unity $1_R$ and also that determinant functions are always constant with value $1_R$.  For example, an FRF $\calF \co \wh{\Q} \to \C$ with determinant 1 which satisfies the Markov condition
\[ \calF(0)^2+\calF(1/0)^2+ \calF(1)^2= \calF(0) \, \calF(1/0) \, \calF(1)\]
is called a {\bf Markov map} in \cite{Bow}.   As in \cite{Bow}, if $\{ \alpha, \beta, \gamma\}$ is a Farey triple
and $\{ a,b,c\} \subset R$ then there is a unique FRF $\calF \co \wh{\Q} \to R$ with
\begin{align*}
\calF(\alpha)&=a & \calF(\beta)&=b & \calF(\gamma)&=c.
\end{align*}

The next lemma shows that the sequence $\calF(\bound(\alpha))$ is a bi-infinite linear recursive sequence with recursion matrix \[W_\alpha = \begin{pmatrix} 0&1 \\ -1 & \calF(\alpha) \end{pmatrix}\] and the reverse of this sequence is also recursive with the same recursion matrix.

\begin{lemma} Suppose that $\calF$ is an FRF and $\alpha \in \wh{\Q}$.  Let $\{ a_j \}_{-\infty}^\infty$ be the ordered bi-infinite sequence obtained by applying $\calF$ to $\bound(\alpha)$.  Then, for every integer $j$,
\[ a_{j+1} = -a_{j-1}+\calF(\alpha) \, a_{j} \quad \text{and} \quad a_{j-1} = -a_{j+1}+\calF(\alpha) \, a_{j} .\]
\end{lemma}
\proof
First, we assume that the indices for $\{ a_j\}$ are shifted so that $\calF(\kappa_-)=a_0$ and $\calF(\kappa_+)=a_1$, where $\kappa_\pm$ are the corners of $\Delta(\alpha)$.

Notice that the first equation in the statement of the lemma holds if and only if the second does.  So, to prove the lemma, it suffices to verify the equalities
\[a_{1}= -a_{-1}+\calF(\alpha) \, a_0 \quad \text{and} \quad a_{2}= -a_{0}+\calF(\alpha) \, a_1\]
occurring at the corners of $\bound(\alpha)$.

The second term of $\bound_-(\alpha)$ coincides with the term in $\bound_+(\kappa_-)$ immediately following $\kappa_+$ and  $\alpha$.  Similarly, the second term of $\bound_+(\alpha)$ immediately follows $\kappa_-$ and $\alpha$ in the sequence $\bound_-(\kappa_+)$.  Hence, the Farey recursive condition implies that
\[ a_{-1}=-a_1+\calF(\alpha) \, a_0 \quad \text{and} \quad a_2=-a_0+\calF(\alpha) \, a_1.\]
These are evidently equivalent to the equalities above. 
\endproof

Define the {\bf generic FRF} to be the FRF \[\calU \co \wh{\Q} \to \Z[x,y,z]\] determined by  
\begin{align*} \calU (0)&=x & \calU (1/0)&=y & \calU (1)&=z. \end{align*}
Consider a composition of $\calU$ with a ring homomorphism $\Z[x,y,z] \to R$.  Since the homomorphism preserves the Farey recursive conditions, the composition will also be an FRF.
On the other hand, if $\calF \co \wh{\Q} \to R$ is an FRF there is a unique ring homomorphism $\psi \co \Z[x,y,z] \to R$ making the diagram 
\[\begin{tikzcd}[column sep=large]
\wh{\Q} \arrow{r}{\calU } \arrow{rd}[swap]{\calF} 
& \Z[x,y,z] \arrow{d}{\psi}\\
& R
\end{tikzcd}
\]
commute.
Hence, every FRF can be viewed as a specialization of $\calU$.

Recall, from Subsection \ref{ssec: FSSBD}, the action of the group $G$ on $\wh{Q}$.  Since an element $\phi \in G$ preserves Farey pairs, the map $\calU  \circ \phi$ is a UFR.  Hence, we have a homomorphism $\psi_\phi \co \Z[x,y,z] \to \Z[x,y,z]$ and a commutative diagram

\[\begin{tikzcd}[column sep=large]
\wh{\Q} \arrow{r}{\calU} \arrow{d}[swap]{\phi} 
& \Z[x,y,z] \arrow{d}{\psi_\phi}\\
\wh{\Q} \arrow{r}{\calU} 
& \Z[x,y,z]
\end{tikzcd}
\]

The homomorphism $\psi_\phi$ is given by evaluating $x,y,z$  at the polynomials 
\[\calU(\phi \cdot 0), \quad \calU(\phi \cdot (1/0)), \quad  \calU(\phi \cdot 1).\]  Because $\phi$ is a bijection, $\psi_\phi$ is an isomorphism with inverse $\psi_{\phi^{-1}}$.

In the cases of 
\[\sigma= \begin{bmatrix} -1&1 \\ 0 & 1 \end{bmatrix}  \quad \text{and} \quad \zeta=  \begin{bmatrix} 0&1 \\ -1 & 1 \end{bmatrix},\] 
$\calU (\gamma)$ and $\calU (\sigma \cdot \gamma)$ differ by interchanging $x$ and $z$, while $\calU (\gamma)$ and $\calU (\zeta \cdot \gamma)$ are related by the cyclic permutation $(x \ z \ y)$ of the three variables.  The $\langle \sigma, \zeta \rangle$-orbit of any $\alpha \in \wh{Q}$ intersects $\Q_0$, which is part of the reason for our focus on this portion of $\calG$.

As the denominator of $\alpha$ grows, the polynomial $\calU(\alpha)$ gets complicated.  The next result is a simple application of the relationship between $\calU$ and the action $G \times \wh{\Q} \to \wh{\Q}$ which shows that the isomorphism class of the affine variety determined by $\calU(\alpha)$ is constant.

\begin{prop}\label{prop: U is irreducible}
If $p \in \calU(\wh{\Q})$ then $p$ is irreducible over $\C$.  More specifically, the affine algebraic variety determined by $p$ is isomorphic to $\C^2$.
\end{prop}
\proof
Suppose $p=\calU(\alpha)$ for some $\alpha \in \wh{\Q}$.  Then there exists $\phi \in G$ with $\phi \cdot 0=\alpha$.  By definition, the isomorphism $\psi_\phi \co \Z[x,y,z] \to \Z[x,y,z]$ takes $x$ to $\calU(\alpha)$.  So we have isomorphic varieties
\[ \C^2 \cong \{ (x,y,z) \in \C^3 \, | \, x=0 \} \cong \{ (x,y,z) \in \C^3 \, | \, p(x,y,z)=0\}. \]
\endproof

\subsection{Traces for $\mathbf{\Omega(\alpha)}$}  \label{ssec: traces} Suppose that $\rho \co \pie (\T) \to \text{SL}_2 \C$ is a representation with character $\chi_\rho$.  Let
\[ x_0= \chi_\rho(a), \quad y_0= \chi_\rho(B), \quad  \text{and} \quad z_0= \chi_\rho(aB).\]
Define the homomorphism $\psi_\rho \co \Z[x,y,z] \to \C$ to be evaluation at $(x_0,y_0,z_0) \in \C^3$.

Take $\beta \in \bound(\alpha)$.  Using the Cayley-Hamilton theorem, the definition of $\Omega$, and the fact that $\chi_\rho$ is constant on conjugacy classes, we find that
\begin{align*} \label{eq: CH} 
\left(\chi_\rho \circ \Omega\right)\left(\alpha \oplus \alpha \oplus \beta\right) &= -\left(\chi_\rho \circ \Omega\right)\left( \beta \right) + \left(\chi_\rho \circ \Omega\right)\left(\alpha\right) \left(\chi_\rho \circ \Omega\right)\left(\alpha \oplus \beta \right).
\end{align*}
In other words, the trace function $\chi_\rho \circ \Omega$ is a UFR and 
\[ (\chi_\rho \circ \Omega)(\alpha) = \psi_\rho \circ \calU (\alpha).\]
The set of all such FRFs is precisely the set of Markov maps from \cite{Bow}.
 
Any two elements of $\pie(\T)$ which are represented by loops freely homotopic to a closed curve with slope $\alpha$ are conjugage in $\pie(\T)$.  So, if $\omega \in \pie \T$ is such an element, then $\chi_\rho(\omega)=\psi_\rho \circ \calU (\alpha)$.  If instead, $\rho$ is a representation into $\text{PSL}_2 \C$, then $x_0$, $y_0$, and $z_0$ are only defined up to sign, as is $\tr(\rho(w))$.  If we choose these signs arbitrarily, then $ \tr(\rho(\omega)) = \pm (\psi_\rho \circ \calU (\alpha))$.

In this paper, we are concerned with representations of 2-bridge link groups.  As we will see, this leads us to consider representations with $y_0=0$.  Let \[\psi_\calT \co \Z[x,y,z] \to \Z[x,z]\] be the evaluation $y=0$.  Define the FRF
\[\calT \co \wh{\Q} \to \Z[x,z]\]
by $\calT = \psi_\calT \circ \calU$.

\begin{lemma}  \label{lem: total degree}
If $p/q \in \wh{\Q}_0$, the total degree of $\calT\left( p/q \right)$ is $q$.
\end{lemma}
\proof
Let $T$ be the UFR obtained from $\calT$ by setting $z=x$. To prove the lemma, we induct on $q$ to show that the degree of $T\left( p/q \right)$ is $q$.
 
Since $T(0)=T(1)=x$ and $T(1/0)=0$, the lemma holds when $q$ is $0$ or $1$.  So, we assume that  $q\geq2$ and that the lemma holds for elements of $\Q_0$ whose denominator is smaller than $q$. 

By Lemma \ref{lem: properside}, we can find a Farey pair $\{ \alpha, \beta \}$ with $\alpha \in \Q_0$ and $\beta \in \wh{\Q}_0$ so that 
\[ p/q = \alpha \oplus \alpha \oplus \beta \]
and the denominators of $\alpha$, $\beta$, and $\alpha \oplus \beta$ are are all less than $q$.  

By definition of $T$,
\[ T\left(p/q\right) = -T(\beta)+T(\alpha) \, T(\alpha \oplus \beta). \] 
By assumption, the degree of $T(\beta)$ is the same as the denominator of $\beta$ and the degree of $T(\alpha) \, T(\alpha \oplus \beta)$ is $q$.  Since the first of these quantities is smaller than the second, the degree of $T(p/q)$ is $q$ as desired.
\endproof

The function $\calT$ has predictable factors of $x$ and $z$.  If we cancel these factors, the resulting polynomial has exponents which are exclusively even.  The next two lemmas make this precise.  To  begin, define  functions $f_j \co \wh{\Q} \to \{ 0, 1\}$ by
\[ f_1(p/q) \equiv pq+1 \pmod{2} \qquad \text{and} \qquad f_2(p/q) \equiv p \pmod{2}.\]
Both functions preserve Farey sums in the sense that, if $\{ \alpha, \beta \}$ are a Farey pair, then
\begin{align} \label{eq: f_j}  
f_j(\alpha \oplus \beta) &\equiv f_j(\alpha) + f_j(\beta) \pmod{2}.
\end{align}
Define $X=x^2$ and $Z=z^2$.  Also, for $\alpha \in \wh{\Q}$, define $f(\alpha)=x^{f_1(\alpha)} z^{f_2(\alpha)}$. 

\begin{lemma} \label{lem: f quot}
If $\{ \alpha, p/q \}$ is a Farey pair in $\wh{\Q}$, then
\[ f\left(p/q\right) = f\left(\alpha \oplus \alpha \oplus p/q\right) \quad \text{and} \quad \frac{f(\alpha) \ f\left(\alpha \oplus p/q\right)}{f\left(p/q\right)} = \begin{cases} 1 & \text{if } q \text{ is even}\\ X & \text{if } pq \text{ is odd} \\ Z & \text{if } p \text{ is even.}\end{cases}\]
\end{lemma}

\proof
The first assertion is a consequence of Equation \ref{eq: f_j} above.  By considering the different possibilities for the parities of $p$ and $q$, the second follows. 
\endproof

Define the {\bf 2-bridge character function} as 
\[\calT_0(\alpha) = \frac{\calT(\alpha)}{f(\alpha)}.\]
It is easy to check that $\calT_0$ is not a UFR.

\begin{lemma}  \label{lem: fac} 
\TO
\end{lemma}
\proof
Here again, we induct on $q$.  Since $\calT(1/0)=0$, the lemma holds for $q=0$.   The recursion matrix for $\calT$ on $\bound(1/0)$ is $\left( \begin{smallmatrix} 0&1 \\ -1 & 0 \end{smallmatrix} \right)$ and it follows that $\{ \calT(n) \}_{n \in \Z}$ is the repeated concatenation of the sequence $\{x, z, -x, -z\}$ shifted so that $\calT(0)=x$.  Since $f(n)$ is $x$ if $n$ is even and $y$ otherwise, the lemma holds for $q=1$.  So, we assume that  $q\geq2$ and that the lemma holds for elements of $\wh{\Q}$ whose denominator is smaller than $q$. 

As in the previous lemma, Lemma \ref{lem: properside} provides a Farey pair $\{ \alpha, \beta \}$ with $\alpha \in \Q$ and $\beta \in \wh{\Q}$ so that 
\[ p/q = \alpha \oplus \alpha \oplus \beta \]
and the denominators of $\alpha$, $\beta$, and $\alpha \oplus \beta$ are are all less than $q$.  

By definition of $\calT$,
\[ \calT\left(p/q\right) = -\calT(\beta)+\calT(\alpha) \, \calT(\alpha \oplus \beta). \] 
From Lemma \ref{lem: f quot}, $f(p/q)=f(\beta)$.  So, dividing through by this, we have
\[ \calT_0\left(p/q\right) = -\calT_0(\beta) + \frac{f(\alpha) f(\alpha \oplus \beta)}{f(\beta)} \, \calT_0(\alpha) \calT_0(\alpha \oplus \beta).\]
Lemma \ref{lem: f quot} and the inductive assumption imply that this is an element of $\Z[X,Z]$.
\endproof

The next lemma follows from Lemma \ref{lem: total degree} and the definition of $\calT_0$.

\begin{cor}\label{cor: degree T_0}
If $p/q \in \Q_0$, the total degree of $\calT_0 \left( p/q \right)$ is $\lfloor \frac{q-1}{2} \rfloor$.
\end{cor}

\section{2-bridge links} \label{ss: 2-bridge} Every 2-bridge link is determined by a number $\alpha \in \Q_0$.  We denote this link as $L(\alpha)$.  Its complement is $M(\alpha) = S^3-L(\alpha)$.  It is well known that the compact manifold formed by removing an open tubular neighborhood of $L(\alpha)$ from $S^3$ can also be obtained by attaching a pair of 2-handles to a thickened 4-holed sphere along a pair of simple closed curves with slopes $1/0$ and $\alpha$.  If we attach only the slope $\infty$ handle, the result is a handle body $H$.  Define
\[ \Gamma_H= \pie(\S)/\normal{b^2} = \big\langle k_0, k_1, k_2, k_3 \, \mid \, k_2=k_1^{-1}, \ k_3=k_0^{-1} \big\rangle = \langle k_0, k_1 \rangle \]
and
\[\Gamma_\alpha = \big\langle k_0, k_1 \, \big| \, \Omega(\alpha)^2 \big\rangle.\]
Then the inclusions of the 4-holed sphere into $H$ and $M(\alpha)$ induce natural identifications of $\Gamma_H$ and $\Gamma_\alpha$ with the fundamental groups of $H$ and $M(\alpha)$.

\begin{figure}[h] 
   \centering
   \includegraphics[width=2.5in]{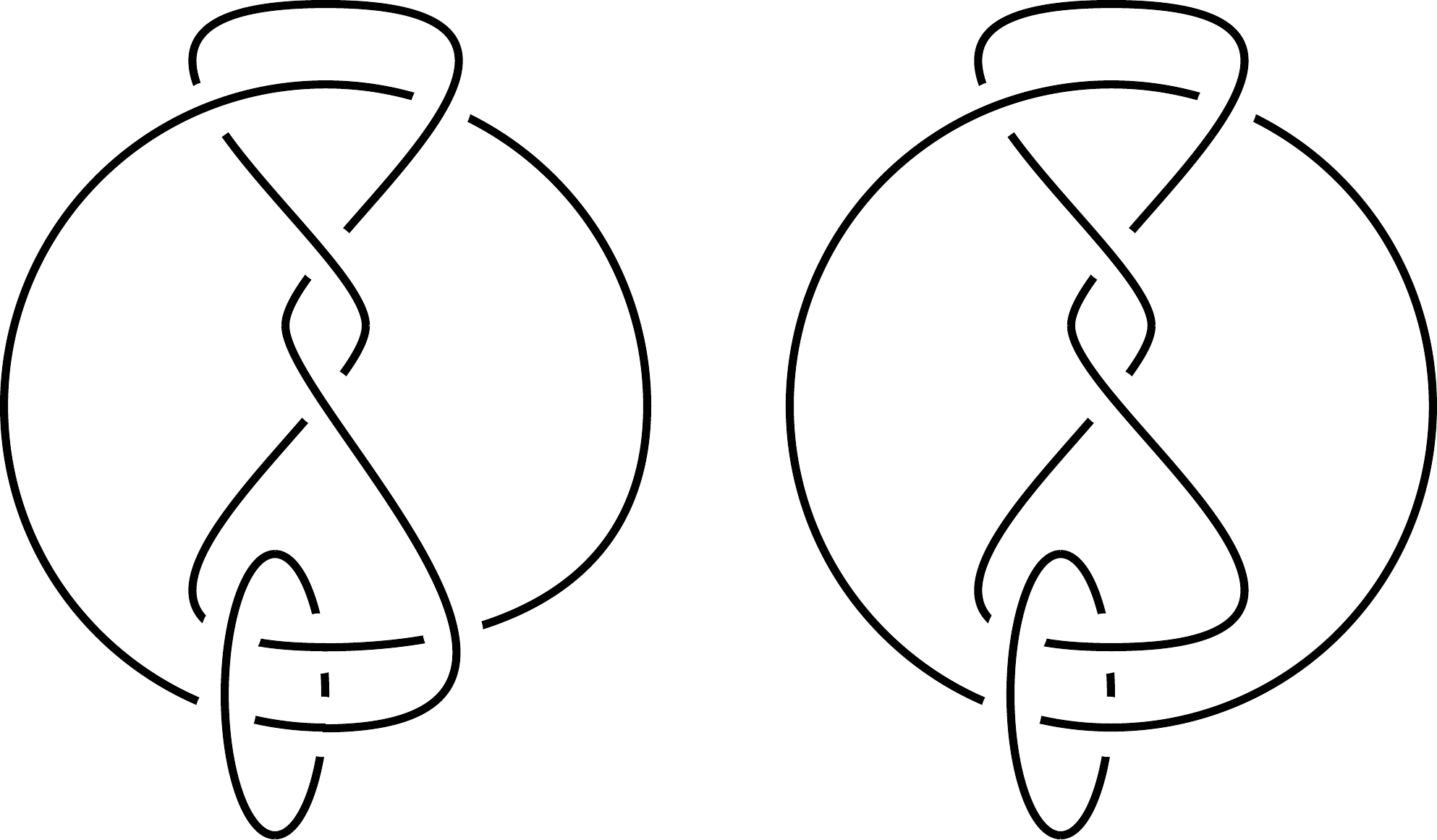} 
   \caption{The manifolds obtained by doing $(1,n)$-Dehn filling on the unknotted cusps are precisely the 2-bridge links which correspond to the terms of $\bound(2/5) \subset \calG$.}   \label{fig: boundary links}
\end{figure}

\begin{remarks} \label{rem: Dehn surgery}
\[ \]
\vspace{-.35in}
\begin{enumerate}
\item It is well-known that $L(\alpha)$ is a knot when the denominator of $\alpha$ is odd and is a 2-component link otherwise.
\item If $\alpha \in \Q_0$ then there is a link complement $M_\alpha$ such that the set of knot complements in $\{ M(\gamma) \, | \, \gamma \in \bound(\alpha) \}$ is the same as the set of $(1,n)$-Dehn fillings on a particular unknotted component (crossing circle) of $M_\alpha$.  Also, by adjusting $M_\alpha$ with a half twist through this crossing circle, we get a similar statement for the 2-component link complements in this set.  See Figure \ref{fig: boundary links}.
\end{enumerate}
\end{remarks}

The next lemma shows that the representations $\rho_{{\bf z},t}$ from Section \ref{sec: reps} are relevant to our discussion of $\Gamma_H$ and $\Gamma_\alpha$.

\begin{lemma} \label{lem: descent}
For every parameter pair $({\bf z},t)$, the representation $\rho_{{\bf z},t}$ descends to $\Gamma_H$.  Moreover, if ${\bf z}$ satisfies $\calT(\alpha)$ then $\rho_{{\bf z},t}$ descends further to $\Gamma_\alpha$.
\end{lemma}

\proof
$\tr\left( \rho_{{\bf z},t} (b)\right)=0$, so $\rho_{{\bf z},t} (b)$ has order two and $\rho_{{\bf z},t}$ descends to $\Gamma_H$.  If ${\bf z}$ satisfies $\calT(\alpha)$ then $\tr\left( \rho_{{\bf z},t} \left( \Omega(\alpha)\right)\right)=0$ and $\rho_{{\bf z},t}$ descends further to $\Gamma_\alpha$.
\endproof

\section{Character varieties}
We refer readers to \cite{CS} for background on $\text{SL}_2 \C$ character varieties and to \cite{BZ}, \cite{GM}, or \cite{LR} for background on $\text{PSL}_2 \C$ character varieties.

Let $\sfX_H$ and $\sfX_\alpha$ be the affine algebraic sets of $\text{SL}_2 \C$ characters of $\Gamma_H$ and $\Gamma_\alpha$ whose algebraic components contain irreducible characters.  Because $\Gamma_\alpha$ is a quotient of $\Gamma_H$, we may regard $\sfX_\alpha$ as a subset of $\sfX_H$.  

As in \cite{CS} and  \cite{GM}, it is typical to use the functions
\[ u= \chi(k_0), \quad w=\chi(k_1), \quad \text{and} \quad  v=\chi(k_0k_1) \]
as coordinates for $\chi$ in $\sfX_H$.   In this paper, we are interested in the {\bf diagonal subvarieties} $\sfD_H \subseteq \sfX_H$ and $\sfD_\alpha \subseteq \sfX_\alpha$ cut out by the polynomial $u-w$.

\begin{remarks} \label{rem: knots}
\[ \]
\vspace{-.35in}
\begin{enumerate}
\item If $L(\alpha)$ is a knot then $k_0$ and $k_1$ are conjugate in $\Gamma_\alpha$ and $\sfD_\alpha=\sfX_\alpha$.  
\item If $L(\alpha)$ is a link, then every component of $\sfX_\alpha$ has dimension two.  Since the components of $\sfD_\alpha$ have dimension one, $\sfD_\alpha$ is a proper subset of $\sfX_\alpha$.  
\item Whenever $L(\alpha)$ is hyperbolic, $\sfD_\alpha$ contains its holonomy characters.
\end{enumerate}
\end{remarks}

The projection $(u,w,v) \to (u,v)$ from $\sfD_H$ to $\C^2$ is a bijection onto its image.  So, by using $(u,v)$ as coordinates, we regard $\sfD_H$ and $\sfD_\alpha$ as affine subsets of $\C^2$.  

\begin{lemma} \label{lem: D(H)}
$\sfD_H = \C^2$.  Moreover, a character $(u_0,v_0) \in \sfD_H$ is reducible if and only if
$(2-v_0)(2+v_0-u_0^2)=0$.
\end{lemma}
\proof
Let $(u_0,v_0) \in \C^2$, take ${\bf z}=(x,z)$ so that 
\begin{align*}
x^2&=2-v_0 &z^2&=2+v_0-u_0^2,
\end{align*}
and take $t \in \C$ to satisfy Equation \ref{eq: Markov} from Section \ref{sec: reps}.  Lemma \ref{lem: descent} implies that $\rho_{{\bf z}, t} \co \Gamma_H \to \text{PSL}_2 \C$ is a representation.  Using Equation \ref{eq: Markov}, it is possible to choose lifts of $\rho_{{\bf z},t}(k_0)$ and $\rho_{{\bf z},t}(k_1)$ to obtain a lift of $\rho_{{\bf z},t}$ to $\text{SL}_2 \C$ whose character $\chi$ satisfies $\chi(k_0)=u_0$ and $\chi(k_0k_1)=v_0$.  The second part of the lemma holds by definition of $\rho_{{\bf z},t}$ (see the last sentence in Section \ref{sec: reps}).
\endproof

Let $\wh{\calT}_0(\alpha) \in \Z[u,v]$ be the polynomial obtained from $\calT_0(\alpha)$ by setting $X=2-v$ and $Z=2+v-u^2$.

\begin{thm} \label{thm: SL}
A character $\chi =(u_0, v_0) \in \C^2$ is an irreducible character in $\sfD_\alpha$ if and only if $\chi$ satisfies $\wh{\calT}_0(\alpha)$ but not $(2-v)  (2+v-u^2)$.
\end{thm}

\proof
By Lemma \ref{lem: D(H)},  $\chi$ is reducible if and only if it satisfies $(2-v)  (2+v-u^2)$.  Together with Lemma \ref{lem: descent}, the proof of Lemma \ref{lem: D(H)} shows that if $\chi$ satisfies $\wh{\calT}_0(\alpha)$ then $\chi \in \sfD_\alpha$. 

It remains only to show that every irreducible character in $\sfD_\alpha$ satisfies $\wh{\calT}(\alpha)$.  Assume, for a contradiction, that this is not the case.  Using Proposition 3.2.1 of \cite{CS}, there is an irreducible plane curve $\sfC \subset \sfD_\alpha$ on which $\wh{\calT}_0(\alpha)$ and $(2-v)(2+v-u^2)$ each have at most finitely many zeros.  Suppose that $(u_0, v_0) \in \C$ avoids these zeros.

As in the proof of Lemma \ref{lem: D(H)}, take ${\bf z}=(x_0,z_0)$ so that 
\begin{align*}
x_0^2&=2-v_0 &z_0^2&=2+v_0-u_0^2,
\end{align*}
and take $t \in \C$ to satisfy Equation \ref{eq: Markov} with $x_0$ and $z_0$.  As before, the restriction of $\rho_{{\bf z}, t}$ to $\Gamma_H$ lifts to $\text{SL}_2 \C$.  This lift represents the character $(u_0, v_0)$ and so $b^2$ and $\Omega(\alpha)^2$ must map trivially under $\rho_{{\bf z},t}$.  But because, $\wh{\calT}_0(\alpha) \neq 0$, the image of $\Omega(\alpha)$ cannot have order 2.  Therefore, $\Omega(\alpha) \in \ker \rho_{{\bf z}, t}$.  

Since $b^2$ and $\Omega(\alpha)$ map trivially under $\rho_{{\bf z},t}$, the image of $\pie \T$ must be abelian.  The expression $(2-v_0)(2+v_0-u_0^2)$ is not zero and commutator $[a,b]$ is $k_0^2$, so the definition of $\rho_{{\bf z}, t}$ implies that $t=\pm i$.  By Equation \ref{eq: Markov}, we must have $u_0=0$.  Hence, $\sfC$ must be the plane curve $\{ (0,v) \, | \, v \in \C \}$.

We've seen that $\Omega(\alpha) \in \ker \rho_{{\bf z},t}$,   Thus,  $\calT(x_0,z_0)^2=4$.  This implies that, after substituting $x^2=2-v$ and  $z^2=2+v$ into $f(\alpha)^2$, the function $f(\alpha)^2 \, \wh{\calT}_0(\alpha)^2$ is never zero on $\sfC$.  But $f(\alpha)^2 \in \{2-v, 2+v, 4-v^2\}$, so $f(\alpha)^2$ has at least one zero on $\sfC$.  This is a contradiction. 
\endproof

If $\rho \co \Gamma_\alpha \to \text{PSL}_2 \C$ is a representation, the traces of $\rho(k_0)$ and $\rho(k_1)$ are only defined up to sign.  On the other hand, $\tr\left(\rho(k_0) \right)^2$ and $\tr\left( \rho(k_0k_1) \right)$ are well-defined.  In fact, the polynomial map $(u^2, v)$ takes characters in $\sfD_\alpha$ to their induced characters in the $\text{PSL}_2 \C$ character variety $\overline{\sfX}_\alpha$ for $\Gamma_\alpha$.  Let $\overline{\sfD}_\alpha$ be the image of $\sfD_\alpha$ under this map.  From \cite{GM}, we know that every $\text{PSL}_2 \C$ representation for $\Gamma_\alpha$ lifts to $\text{SL}_2 \C$.  This implies that, as in Remark \ref{rem: knots}, if $L(\alpha)$ is a knot, then $\overline{\sfD}_\alpha$ is equal to the entire $\text{PSL}_2 \C$ character variety for $\Gamma_\alpha$.

It is natural to use $X=2-v$ and $Z=2+v-u^2$ to change coordinates on $\overline{\sfD}_\alpha$ and we do so for the remainder of this paper.   The next result follows immediately from Theorem \ref{thm: SL}.

\begin{thm} \label{thm: PSL}
\PSL 
\end{thm}

\subsection{Common factors and epimorphisms} \label{subsec: epi}

Define $G_\alpha$ to be the stabilizer in $G$ of $\alpha\in \wh{\Q}$ and define $\wh{G}_\alpha$ to be the subgroup in $G$ generated by $G_\infty$ and $G_\alpha$.  Let $\wh{G}_\alpha \{ \alpha, \infty\}$ be the union of the $\wh{G}_\alpha$-orbits of $\alpha$ and $\infty$.  The main result of \cite{ORS}, states that if $\alpha' \in \wh{G}_\alpha \{ \alpha, \infty\}$ then there is an epimorphism $\Gamma_\alpha \to \Gamma_{\alpha'}$ which takes $k_j\in \Gamma_\alpha$ to $k_j \in \Gamma_{\alpha'}$.  As an application of Corollary \ref{thm: PSL} and the action $G \times \wh{\Q} \to \wh{\Q}$, we give an elementary argument reproducing this result.

For $\phi \in G$ let $\theta_\phi \co \Z[x,z] \to \Z[x,z]$ be the ring homomorphism determined by \[ \theta_\phi(x) = \calT(\phi \cdot 0)\quad \text{and} \quad \theta_\phi(z) = \calT(\phi \cdot 1).\]  For a polynomial $p \in \C[x,z]$, let $\pi_p \co \C[x,z] \to \C[x,z]/(p)$ be the natural quotient homomorphism. The next pair of results describe the relationship between $\calT$ and the action $G \times \wh{\Q} \to \wh{\Q}$.

\begin{lemma} \label{lem: T action}
If $\phi \in G$, $p \in \C[x,z]$, and $\calT(\phi \cdot (1/0)) \in \ker \pi_p$ then 
\[ \pi_p \circ \calT \circ \phi = \pi_p \circ \theta_\phi \circ \calT. \]
\end{lemma}

\proof
The map $\pi_p \circ \calT \circ \phi$ is a UFR, so there is a unique ring homomorphism $\psi \co \Z[x,y,z] \to \C[x,z]/(p)$ such that $\pi_p \circ \calT \circ \phi = \psi \circ \calU$.  The map $\psi$ is uniquely determined by 
\[ \psi(x)=\pi_p \circ \calT(\phi \cdot 0), \quad \psi(y)=\pi_p \circ \calT(\phi \cdot \infty)=0, \quad \text{and} \quad \psi(z)=\pi_p \circ \calT(\phi \cdot \infty).\]
Hence $\psi = \pi_p \circ \theta_\phi \circ \psi_\calT$.  Since $\calT= \psi_\calT \circ \calU$, this proves the lemma.
\endproof

\begin{prop}\label{prop: factors}
If $\alpha \in \wh{\Q}$ and $p \in \C[x,z]$ such that $\calT(\alpha) \in \ker \pi_p$ then
\[ \pi_p \circ \calT \circ \phi = \pm \pi_p \circ \calT\]
for every $\phi \in \wh{G}_\alpha$.
\end{prop}
\proof
To prove the theorem, it is enough to show that the conclusion holds separately for $\phi \in G_\infty$ and for $\phi \in G_\alpha$.

Let 
\[ \phi_1 = \begin{bmatrix} -1 & 0 \\ 0 & 1 \end{bmatrix}  \quad \text{and} \quad \phi_2= \begin{bmatrix} -1 & 2 \\ 0 & 1 \end{bmatrix}\]
and notice that $\phi_1$ and $\phi_2$ generate $G_\infty$.  Also,
\begin{align*}
\phi_1 \cdot 0 & = 0 & \phi_2 \cdot 0 &= 2\\
\phi_1 \cdot 1 & = -1 & \phi_2 \cdot 1&=1.
\end{align*}
By the second sentence in the proof of Lemma \ref{lem: fac}, $\theta_{\phi_1}$ fixes $x$ and reverses the sign of $z$.  Similarly, $\theta_{\phi_2}$ fixes $z$ and reverses the sign of $z$.  For $j=1,2$, $\calT (\phi_j \cdot \infty) = 0$, so Lemmas \ref{lem: fac} and \ref{lem: T action} imply that
\begin{align} \label{eq: G_infty} \pi_q \circ \calT \circ \phi_j &= \pm \pi_q \circ \calT\end{align}
for every $q \in \C[x,z]$.
Since $G_\infty = \langle \phi_1, \phi_2 \rangle$, the same formula holds if we replace $\phi_j$ with any $\phi \in G_\infty$. 

Now take $\psi \in G$ so that $\psi \cdot \alpha = \infty$ and note that $G_\alpha = \psi^{-1} G_\infty \psi$.  Let $\theta$ denote both the isomorphism $\theta_{\psi^{-1}}$ as well as the induced map on quotient rings.
\[\begin{tikzcd}[column sep=large]
\C[x,z] \arrow{r}{\theta} \arrow{d}[swap]{\pi_p} 
& \C[x,z] \arrow{d}{\pi_{\theta(p)}}\\
\C[x,z]/(p) \arrow{r}{\theta} 
& \C[x,z]/(\theta(p))
\end{tikzcd}
\]
Since $\psi^{-1} \cdot \infty = \alpha$, Lemma \ref{lem: T action} and the diagram above give
\begin{align} 
\pi_p \circ \calT \circ \psi^{-1}  &= \pi_p \circ \theta \circ \calT \nonumber \\
&= \theta \circ \pi_{\theta(p)} \circ \calT. \label{eq: psi}
\end{align}
In particular,
\begin{align*} 
\pi_p \circ \calT &= \pi_p \circ \calT \circ \psi^{-1} \circ \psi \\
&= \theta \circ \pi_{\theta(p)} \circ \calT \circ \psi.
\end{align*}
Take $\xi \in G_\alpha$.  Then $\xi = \psi^{-1} \phi \psi$ for some $\phi \in G_\infty$.  Then,  using Equations \ref{eq: G_infty} and \ref{eq: psi}, we have
\begin{align*}
 \pi_p \circ \calT \circ \xi &= \theta \circ \pi_{\theta(p)} \circ\calT \circ \phi  \circ \psi \\
 &= \pm \theta \circ \pi_{\theta(p)} \circ \calT \circ \psi \\
 &= \pm \pi_p \circ \calT.
 \end{align*}
\endproof

The next corollary is immediate.

\begin{cor}\label{cor: factors}
\Factors 
\end{cor}

\begin{cor} \label{cor: epi}  
\Epimorphism 
\end{cor}

\proof
Since $M(\alpha)$ is hyperbolic, we may take $\chi =(X_0, Z_0) \in \overline{\sfD}_\alpha$ to be the character of a holonomy representation for $\Gamma_\alpha$.  By Corollary \ref{thm: PSL}, $\chi$ satisfies $\calT_0(\alpha)$ and $X_0Z_0 \neq 0$.  Proposition \ref{prop: factors} implies that $\chi$ also satisfies $\calT_0(\alpha')$.  Corollary \ref{thm: PSL} implies that $\chi$ also corresponds to a character for a representation for $\Gamma_{\alpha'}$.  

Since $\chi$ is an irreducible character for $\Gamma_H$, we may choose a single representation $\rho \co \Gamma_H \to \text{PSL}_2 \C$ which represents $\chi$ and factors through both $\Gamma_\alpha$ and $\Gamma_{\alpha'}$.  The corresponding representation $\Gamma_\alpha \to \text{PSL}_2 \C$ is faithful, so we obtain an epimorphism as claimed in the statement of the corollary.
\endproof

\begin{quest}
Lee and Sakuma show in \cite{LS} that the converse of Corollary 7.7 holds.  Is there an elementary proof of this using Farey recursion?
\end{quest}

\subsection{Relationship to \cite{MPV}} \label{subsec: MPV}
Our work here can be easily reconciled with that in \cite{MPV}, wherein the authors discuss the the character varieties of double twist knots.  They describe these knots with a pair of integers $k$ and $l$, where $l$ is even and $k>0$, and denote them as $J(k,l)$.  The following facts can be shown with straightforward induction arguments.

\begin{enumerate}
\item The group presentation used by \cite{MPV} for the knot $J(k,l)$ corresponds exactly to our presentation for $\Gamma_{\gamma_{k,l}}$ where 
\[ \gamma_{k,l} = \frac{1+l (k-1)}{1+lk}.\]
For fixed $k$, these numbers form the sequence obtained from $\bound\left(\frac{k-1}{k}\right)$ by deleting every other term in such a way that the corner $v_1$ survives.
\item The polynomial given in Proposition 3.8 of \cite{MPV} which determines the $\text{PSL}_2 \C$ character variety for $J(k,l)$ becomes $\pm \calT_0(\gamma_{k,l})$ after letting $y=2-X-Z$ and $r=2-Z$.
\end{enumerate}

In view of the recursion matrices for $\calT$, the substitution $s=\calT(\alpha)$ is natural while considering the terms of $\bound(\alpha)$.  Lemma \ref{lem: fac} tells us that the value $S=s^2$ can be expressed as a polynomial in $\Z[X,Z]$.  For $J(k,l)$, we have $S=\calT\left( \frac{k-1}{k} \right)^2$.  Since the numbers $\frac{k-1}{k}$ are in $\bound(1)$, their images under $\calT$ satisfy the recursion given by $\left( \begin{smallmatrix} 0 & 1 \\ -1 & z \end{smallmatrix}\right)$ and induction shows that $S=\calT\left( \frac{k-1}{k} \right)^2$ is linear in $X$.  So, as shown in \cite{MPV}, the map $(X,Z) \mapsto (S,Z)$ gives a birational change of coordinates for $\overline{\sfD}_{\gamma_{k,l}}$.
\begin{enumerate}
\item[(3)] In Section 4 of \cite{MPV}, the authors introduce polynomials in $\Z[r,t]$ whose corresponding plane curves have projective closures in $\mathbb{P}^1 \times \mathbb{P}^1$ which realize the smooth models for the double twist knots.  If we set $t=2-S$ and $r=2-Z$ in their polynomials, the result is the same (up to sign) as when we change $X$ to $S$ in $\calT_0\left( \gamma_{k,l} \right)$ as described above.
\end{enumerate}
Item (3) suggests that it might be fruitful to consider the maps $(X,Z) \mapsto (S,Z)$ and the techniques of \cite{MPV} more generally across the 2-bridge link family.

\section{Riley polynomials and holonomy representations} \label{sec: Riley}
Suppose $\alpha \in \Q_0$.  A {\bf p-rep} for $\Gamma_\alpha$ is an irreducible $\text{PSL}_2 \C$ representation for $\Gamma_\alpha$ whose character satisfies $u^2=w^2=4$ (or equivalently $Z=-X$).   In particular, characters of p-reps for $\Gamma_\alpha$ lie on $\overline{\sfD}_\alpha$.  

Define $\psi_\Lambda \co \Z[X,Z] \to \Z[X]$ to be the ring homomorphism defined by $X \mapsto X$ and $Z \mapsto -X$. Define 
\[\Lambda_\alpha = \psi_\Lambda \circ \calT_0(\alpha).\]
The next result is a direct consequence of Corollary \ref{thm: PSL}.

\begin{cor} \label{cor: p-reps}
A point $(X_0,-X_0) \in \C^2$ is the character of of an irreducible p-rep for $\Gamma_\alpha$ if and only if $X_0 \neq 0$ and $\Lambda_\alpha(X_0)=0$.  When this is the case, $(X_0,-X_0)$ is a character for the representation determined by \begin{align*}
k_0 &\mapsto \begin{bmatrix} 1 & 1 \\ 0 & 1 \end{bmatrix} & k_1 &\mapsto \begin{bmatrix} 1 & 0 \\ -X_0 & 1 \end{bmatrix}. 
\end{align*}
\end{cor}

The {\bf Riley polynomial} for a number $\alpha \in \Q_0$ with odd denominator is defined in \cite{RR}.  In Theorem 2 of \cite{RR}, Riley shows that the Riley polynomial for $\alpha$ is monic and has degree $\frac{1}{2}(q-1)$, where $q$ is the denominator of $\alpha$.  Also, as part of this this theorem, Riley shows that $X_0$ satisfies the Riley polynomial for $\alpha$ if and only if the assignment
\begin{align*}
k_0 &\mapsto \begin{bmatrix} 1 & 1 \\ 0 & 1 \end{bmatrix} & k_1 &\mapsto \begin{bmatrix} 1 & 0 \\ -X_0 & 1 \end{bmatrix}
\end{align*}
 determines a p-rep for $\Gamma_\alpha$.  
  
 Corollary \ref{cor: degree T_0} shows that the degree of $\Lambda_\alpha$ is the same as the degree of the Riley polynomial, which establishes the next result.  
 
 \begin{cor} \label{cor: Riley polynomial}
 \Riley
  \end{cor}
 
Suppose $\alpha=\frac{p}{q} \in \Q_0$. It follows from Thurston's geometrization theorems that $M(\alpha)$ admits a complete hyperbolic structure with finite volume if and only if $p \notin \{1, q-1\}$.  The corresponding holonomy representations for $\Gamma_\alpha$ are always p-reps.  In particular, there is a root of $\calT_0(\alpha)$ which provides the holonomy as in Corollary \ref{cor: p-reps}.   It is well-known that the traces of the holonomy for a knot complement $M(\alpha)$ are algebraic integers.  In fact, whenever the holonomy representation of a hyperbolic 3-manifold has non-integer traces, the manifold contains a closed essential surface.  Hatcher and Oertel show in \cite{HO} that no 2-bridge link contains a closed essential surface and so their holonomy representations are integral.  According to Alan Reid, Riley was also aware that the p-rep polynomials for two component 2-bridge links are monic.
 
 In the present setting, it is not difficult to reproduce this result.  Here, we show that $\Lambda_\alpha$ is monic for every $\alpha \in \Q_0$. 
 
 \begin{thm} \label{thm: monic}
 If $\alpha \in \Q_0$ then $\Lambda_{\alpha}$ is monic.
 \end{thm}
 
 \proof
 Let $d$ be the denominator of $\alpha$.  Since $\Lambda_0=\Lambda_1=\Lambda_{1/2}$ the theorem holds when $d\leq 2$.  Suppose then, that $d\geq 3$ and that the theorem holds for elements of $\Q_0$ whose denominators are smaller than $d$.
 
As usual, Lemma \ref{lem: properside} provides a Farey pair $\{ \gamma, \beta \}$ with $\gamma \in \Q_0$ and $\beta \in \wh{\Q}_0$ so that
\[ \alpha  = \gamma \oplus \gamma \oplus \beta\]
and the denominators of $\gamma$, $\beta$, and $\gamma \oplus \beta$ are all less than $d$.  Also, because $d \geq 3$, $\beta \neq \frac{1}{0}$.  Write $\gamma=p/q$ and $\beta=r/s$.  Define \[ Q= \begin{cases} 1& \text{if } s \text{ is even} \\
X & \text{if } rs \text{ is odd} \\
-X & \text{if } r \text{ is even}.
\end{cases} \]

From Lemma \ref{lem: f quot} and the proof of Lemma \ref{lem: fac}, 
\[ \Lambda_\alpha = -\Lambda_\beta + Q \, \Lambda_\gamma \, \Lambda_{\gamma \oplus \beta}.\]
Using Lemma \ref{cor: degree T_0} to compute degrees, we have
\[ \deg \Lambda_\alpha - \deg \Lambda_\beta = \left\lfloor \frac{2q+s-1}{2} \right\rfloor - \left\lfloor \frac{s-1}{2} \right\rfloor = q\] which is positive.  Therefore, the leading coefficient of $\Lambda_\alpha$ is the same as that of $Q \, \Lambda_\gamma \, \Lambda_{\gamma \oplus \beta}$.  This coefficient is one, by the inductive assumption.
\endproof
 
 \begin{cor}
 If $G \subset \text{PSL}_2 \C$ is a Kleinian group which uniformizes a hyperbolic 2-bridge link then the traces (defined up to sign) of the elements of $G$ are algebraic integers.
 \end{cor}
 
 \proof
 $G$ is conjugate to the group 
 \[ \left\langle  \begin{bmatrix} 1& 1\\ 0 & 1 \end{bmatrix}, \ \begin{bmatrix} 1& 0\\ -X_0 & 1 \end{bmatrix} \right\rangle\]
 where $X_0$ is a non-zero root of $\Lambda_\alpha$ for some $\alpha \in \Q_0$.  The trace of the product of these two matrices is $2-X_0$ and it follows that the traces of $G$ lie in the ring $\Z[X_0]$.  By Theorem \ref{thm: monic}, we know that this ring is made up entirely of algebraic integers.
 \endproof
 
 We conclude with three observations made from the data we've collected.
 
 First, our computer experiments show that, except for $\Lambda_{1/4}$, $\Lambda_{p/q}$ has a non-trivial factorization over $\Q$ for every even $q$ which is less than 300.  In contrast, it is not difficult to find $\alpha \in \Q$ for which $\Lambda_\alpha$ is irreducible.  For instance, Hoste  and Shanahan prove  in \cite{HS3} that this is the case for each of the twist knots $L\left( \frac{n}{2n+1}\right)$.
 
In a similar vein, Theorem 3 of \cite{RR} states that Riley polynomials have no repeated factors.  On the other hand, our experiments find many $\alpha$'s with even denominator less than 300 for which $\Lambda_\alpha$ has factors with multiplicity larger than one.  In almost every such case, the factor with high multiplicity is the monomial $X$.  Amongst the numbers we searched, we found exactly three exceptions to this.  If $\alpha= 7/24$ then $\Lambda_\alpha$ is divisible by $(X^2-1)^2$ and if $\alpha \in  \{41/264, \,103/264\}$ then $\Lambda_\alpha$ is divisible by $(X^2-1)^3$.  We remark that, in \cite{RW}, Theorem 3 of \cite{RR} is used to show that there is no 2-bridge knot whose trace field contains $\Q(i)$ as a subfield.

Finally, Question 1.7 of \cite{MPV} asks whether there is an $\alpha \in \Q$ with odd denominator and numerator not equal to one for which $\calT_0(\alpha)$ is irreducible over $\Q$ but $\Lambda_\alpha$ factors non-trivially.  They found no such examples amongst the double twist knots with crossing numbers at most 98.  Since there are many $\alpha$'s with even denominator (links) for which $\calT_0(\alpha)$ is irreducible, the observation above provides many examples with {\em even} denominators.  However, if we look only at $\alpha$'s with odd denominators (knots), we find no examples up to denominator 200.
 
\section{Lists} \label{sec: lists}

This section consists of lists of polynomials computed using the techniques of this paper.

\begin{center} $T_0(\alpha)$ with $\alpha \in \Q \cap \left[ 0, \frac{1}{2} \right]$ and denominator at most $18$. \end{center} 
\begin{align*} 
1/3 && X - 1 \\
1/4 && X - 2\\
1/5 && X^{2} - 3  X + 1\\
2/5 && X Z - Z - 1\\
1/6 && {\left(X - 1\right)} {\left(X - 3\right)}\\
1/7 && X^{3} - 5  X^{2} + 6  X - 1\\
2/7 && X^{2} Z - 3  X Z + 2  Z - 1\\
3/7 && X^{2} Z - X Z - 2  X + 1\\
1/8 && {\left(X^{2} - 4  X + 2\right)} {\left(X - 2\right)}\\
3/8 && X^{2} Z - 2  X Z - X + Z\\
1/9 && {\left(X^{3} - 6  X^{2} + 9  X - 1\right)} {\left(X - 1\right)}\\
2/9 && X^{3} Z - 5  X^{2} Z + 7  X Z - 2  Z - 1\\
4/9 && X^{2} Z^{2} - X Z^{2} - 3  X Z + 2  Z + 1\\
1/10 && (X^2 - 3X + 1)(X^2 - 5X + 5) \\
3/10 && X^3Z - 4X^2Z + 5XZ - 2X - 2Z + 3 \\
1/11 &&  X^5 - 9X^4 + 28X^3 - 35X^2 + 15X - 1 \\
2/11 &&  X^4Z - 7X^3Z + 16X^2Z - 13XZ + 3Z - 1 \\
3/11 &&  X^4Z - 5X^3Z + 8X^2Z - X^2 - 4XZ + X + 1 \\
4/11 &&  X^3Z^2 - 3X^2Z^2 - X^2Z + 3XZ^2 - Z^2 + Z + 1 \\
5/11 &&  X^3Z^2 - X^2Z^2 - 4X^2Z + 3XZ + 3X - 1 \\
\end{align*}
\begin{align*}
1/12 &&  \scriptstyle (X^2 - 4X + 1)(X - 1)(X - 2)(X - 3) \\
5/12 &&  \scriptstyle (XZ - Z - 2)(XZ - 1)(X - 1) \\
1/13 &&  \scriptstyle X^6 - 11X^5 + 45X^4 - 84X^3 + 70X^2 - 21X + 1 \\
2/13 &&  \scriptstyle X^5Z - 9X^4Z + 29X^3Z - 40X^2Z + 22XZ - 3Z - 1 \\
3/13 &&  \scriptstyle X^5Z - 7X^4Z + 17X^3Z - 16X^2Z - 2X^2 + 4XZ + 5X - 1 \\
4/13 &&  \scriptstyle X^4Z^2 - 5X^3Z^2 + 9X^2Z^2 - 3X^2Z - 7XZ^2 + 8XZ + 2Z^2 - 5Z + 1 \\
5/13 &&  \scriptstyle X^4Z^2 - 3X^3Z^2 - 2X^3Z + 3X^2Z^2 + 3X^2Z - XZ^2 + X^2 - XZ - X + 1 \\
6/13 &&  \scriptstyle X^3Z^3 - X^2Z^3 - 5X^2Z^2 + 4XZ^2 + 6XZ - 3Z - 1 \\
1/14 &&  \scriptstyle (X^3 - 5X^2 + 6X - 1)(X^3 - 7X^2 + 14X - 7) \\
3/14 && \scriptstyle X^5Z - 8X^4Z + 23X^3Z - 28X^2Z - X^2 + 13XZ + 2X - 2Z + 1 \\
5/14 &&  \scriptstyle X^4Z^2 - 4X^3Z^2 - X^3Z + 6X^2Z^2 - 4XZ^2 + 3XZ + Z^2 + 2X - 2Z - 1 \\
1/15 && \scriptstyle (X^4 - 9X^3 + 26X^2 - 24X + 1)(X^2 - 3X + 1)(X - 1) \\
2/15 && \scriptstyle X^6Z - 11X^5Z + 46X^4Z - 91X^3Z + 86X^2Z - 34XZ + 4Z - 1 \\
4/15 && \scriptstyle (X^3Z - 3X^2Z + 2XZ - 1)(X^2Z - 4XZ + 4Z - 1) \\
7/15 && \scriptstyle X^4Z^3 - X^3Z^3 - 6X^3Z^2 + 5X^2Z^2 + 10X^2Z - 6XZ - 4X + 1 \\
1/16 && \scriptstyle (X^4 - 8X^3 + 20X^2 - 16X + 2)(X^2 - 4X + 2)(X - 2) \\
3/16 && \scriptstyle X^6Z - 10X^5Z + 38X^4Z - 68X^3Z + 58X^2Z - 2X^2 - 22XZ + 7X + 3Z - 4 \\
5/16 && \scriptstyle X^5Z^2 - 6X^4Z^2 + 14X^3Z^2 - 4X^3Z - 16X^2Z^2 + 15X^2Z + 9XZ^2 - 18XZ - 2Z^2 + 3X + 7Z - 4 \\
7/16 && \scriptstyle (X^3Z^2 - 2X^2Z^2 - 3X^2Z + XZ^2 + 4XZ + X - Z)(XZ - 2) \\
1/17 && \scriptstyle X^8 - 15X^7 + 91X^6 - 286X^5 + 495X^4 - 462X^3 + 210X^2 - 36X + 1 \\
2/17 && \scriptstyle X^7Z - 13X^6Z + 67X^5Z - 174X^4Z + 239X^3Z - 166X^2Z + 50XZ - 4Z - 1 \\
3/17 && \scriptstyle X^7Z - 11X^6Z + 47X^5Z - 98X^4Z + 103X^3Z - X^3 - 51X^2Z + 3X^2 + 9XZ - 1 \\
4/17 && \scriptstyle X^6Z^2 - 9X^5Z^2 + 31X^4Z^2 - 50X^3Z^2 - 3X^3Z + 36X^2Z^2 + 14X^2Z - 8XZ^2 - 18XZ + 4Z + 1 \\
5/17 && \scriptstyle X^6Z^2 - 7X^5Z^2 + 19X^4Z^2 - 3X^4Z - 25X^3Z^2 + 13X^3Z + 16X^2Z^2 - 18X^2Z - 4XZ^2 + 2X^2 + 8XZ - 4X + 1 \\
6/17 && \scriptstyle X^5Z^3 - 5X^4Z^3 - X^4Z^2 + 10X^3Z^3 - 10X^2Z^3 + 6X^2Z^2 + 5XZ^3 + 3X^2Z - 8XZ^2 - Z^3 - 3XZ + 3Z^2 - 1 \\
7/17 && \scriptstyle X^5Z^3 - 3X^4Z^3 - 4X^4Z^2 + 3X^3Z^3 + 9X^3Z^2 - X^2Z^3 + 5X^3Z - 6X^2Z^2 - 9X^2Z + XZ^2 - 2X^2 + 4XZ + 4X - 1 \\
8/17 && \scriptstyle X^4Z^4 - X^3Z^4 - 7X^3Z^3 + 6X^2Z^3 + 15X^2Z^2 - 10XZ^2 - 10XZ + 4Z + 1 \\
1/18 && \scriptstyle (X^3 - 6X^2 + 9X - 1)(X^3 - 6X^2 + 9X - 3)(X - 1)(X - 3) \\
5/18 && \scriptstyle (X^5Z^2 - 7X^4Z^2 + 18X^3Z^2 - 2X^3Z - 20X^2Z^2 + 7X^2Z + 8XZ^2 - 5XZ + X - 2Z - 1)(X - 1) \\
7/18 && \scriptstyle (X^4Z^3 - 3X^3Z^3 - 3X^3Z^2 + 3X^2Z^3 + 5X^2Z^2 - XZ^3 + 3X^2Z - 2XZ^2 - 3XZ - X + 2Z + 1)(X - 1) \\
\end{align*}


\newpage
\begin{center} Riley polynomials $\Lambda_\alpha(X)$ with $\alpha \in \Q \cap \left[ 0, \frac{1}{2} \right]$ and denominator at most $20$. \end{center} 
\begin{align*}
1/3 && X - 1 \\
1/4 && X - 2 \\
1/5 && X^2 - 3 X + 1 \\
2/5 && -X^2 + X - 1 \\
1/6 && (X - 1) (X - 3) \\
1/7 && X^3 - 5 X^2 + 6 X - 1 \\
2/7 && -X^3 + 3 X^2 - 2 X - 1 \\
3/7 && -X^3 + X^2 - 2 X + 1 \\
1/8 && (X^2 - 4 X + 2) (X - 2) \\
3/8 && -(X^2 - 2 X + 2) X \\
1/9 && (X^3 - 6 X^2 + 9 X - 1) (X - 1) \\
2/9 && -X^4 + 5 X^3 - 7 X^2 + 2 X - 1 \\
4/9 && X^4 - X^3 + 3 X^2 - 2 X + 1 \\
1/10 && (X^2 - 3 X + 1) (X^2 - 5 X + 5) \\
3/10 && -(X^2 - X - 1) (X^2 - 3 X + 3) \\
1/11 && X^5 - 9 X^4 + 28 X^3 - 35 X^2 + 15 X - 1 \\
2/11 && -X^5 + 7 X^4 - 16 X^3 + 13 X^2 - 3 X - 1 \\
3/11 && -X^5 + 5 X^4 - 8 X^3 + 3 X^2 + X + 1 \\
4/11 && X^5 - 3 X^4 + 4 X^3 - X^2 - X + 1 \\
5/11 && X^5 - X^4 + 4 X^3 - 3 X^2 + 3 X - 1 \\
1/12 && (X^2 - 4 X + 1) (X - 1) (X - 2) (X - 3) \\
5/12 && (X^2 - X + 2) (X^2 + 1) (X - 1) \\
1/13 && X^6 - 11 X^5 + 45 X^4 - 84 X^3 + 70 X^2 - 21 X + 1 \\
2/13 && -X^6 + 9 X^5 - 29 X^4 + 40 X^3 - 22 X^2 + 3 X - 1 \\
3/13 && -X^6 + 7 X^5 - 17 X^4 + 16 X^3 - 6 X^2 + 5 X - 1 \\
4/13 && X^6 - 5 X^5 + 9 X^4 - 4 X^3 - 6 X^2 + 5 X + 1 \\
5/13 && X^6 - 3 X^5 + 5 X^4 - 4 X^3 + 2 X^2 - X + 1 \\
6/13 && -X^6 + X^5 - 5 X^4 + 4 X^3 - 6 X^2 + 3 X - 1 \\
1/14 && (X^3 - 5 X^2 + 6 X - 1) (X^3 - 7 X^2 + 14 X - 7) \\
3/14 && -(X^3 - 3 X^2 + 2 X - 1) (X^3 - 5 X^2 + 6 X + 1) \\
5/14 && (X^3 - X^2 + 1) (X^3 - 3 X^2 + 4 X - 1) \\
\end{align*}

\begin{align*}
1/15 &&   \scriptstyle (X^4 - 9 X^3 + 26 X^2 - 24 X + 1) (X^2 - 3 X + 1) (X - 1) \\
2/15 &&   \scriptstyle -X^7 + 11 X^6 - 46 X^5 + 91 X^4 - 86 X^3 + 34 X^2 - 4 X - 1 \\
4/15 &&   \scriptstyle (X^4 - 3 X^3 + 2 X^2 + 1) (X^3 - 4 X^2 + 4 X + 1) \\
7/15 &&   \scriptstyle -X^7 + X^6 - 6 X^5 + 5 X^4 - 10 X^3 + 6 X^2 - 4 X + 1 \\
1/16 &&   \scriptstyle (X^4 - 8 X^3 + 20 X^2 - 16 X + 2) (X^2 - 4 X + 2) (X - 2) \\
3/16 &&   \scriptstyle -(X^4 - 6 X^3 + 10 X^2 - 2 X - 2) (X^3 - 4 X^2 + 4 X - 2) \\
5/16 &&   \scriptstyle (X^4 - 4 X^3 + 6 X^2 - 2 X - 2) (X^3 - 2 X^2 + 2) \\
7/16 &&  \scriptstyle -(X^4 - 2 X^3 + 4 X^2 - 4 X + 2) (X^2 + 2) X \\
1/17 &&  \scriptstyle X^8 - 15 X^7 + 91 X^6 - 286 X^5 + 495 X^4 - 462 X^3 + 210 X^2 - 36 X + 1 \\
2/17 && \scriptstyle  -X^8 + 13 X^7 - 67 X^6 + 174 X^5 - 239 X^4 + 166 X^3 - 50 X^2 + 4 X - 1 \\
3/17 &&  \scriptstyle -X^8 + 11 X^7 - 47 X^6 + 98 X^5 - 103 X^4 + 50 X^3 - 6 X^2 - 1 \\
4/17 && \scriptstyle  X^8 - 9 X^7 + 31 X^6 - 50 X^5 + 39 X^4 - 22 X^3 + 18 X^2 - 4 X + 1 \\
5/17 &&  \scriptstyle X^8 - 7 X^7 + 19 X^6 - 22 X^5 + 3 X^4 + 14 X^3 - 6 X^2 - 4 X + 1 \\
6/17 &&  \scriptstyle -X^8 + 5 X^7 - 11 X^6 + 10 X^5 + X^4 - 10 X^3 + 6 X^2 - 1 \\
7/17 &&  \scriptstyle -X^8 + 3 X^7 - 7 X^6 + 10 X^5 - 11 X^4 + 10 X^3 - 6 X^2 + 4 X - 1 \\
8/17 &&  \scriptstyle X^8 - X^7 + 7 X^6 - 6 X^5 + 15 X^4 - 10 X^3 + 10 X^2 - 4 X + 1 \\
1/18 &&  \scriptstyle (X^3 - 6 X^2 + 9 X - 1) (X^3 - 6 X^2 + 9 X - 3) (X - 1) (X - 3) \\
5/18 && \scriptstyle  (X^4 - 3 X^3 + X^2 + 2 X + 1) (X^3 - 4 X^2 + 5 X - 1) (X - 1) \\
7/18 && \scriptstyle  -(X^4 - X^3 + X^2 + 1) (X^3 - 2 X^2 + 3 X - 1) (X - 1) \\
1/19 &&  \scriptstyle X^9 - 17 X^8 + 120 X^7 - 455 X^6 + 1001 X^5 - 1287 X^4 + 924 X^3 - 330 X^2 + 45 X - 1 \\
2/19 &&  \scriptstyle -X^9 + 15 X^8 - 92 X^7 + 297 X^6 - 541 X^5 + 553 X^4 - 296 X^3 + 70 X^2 - 5 X - 1 \\
3/19 && \scriptstyle  -X^9 + 13 X^8 - 68 X^7 + 183 X^6 - 269 X^5 + 211 X^4 - 80 X^3 + 18 X^2 - 9 X + 1 \\
4/19 &&  \scriptstyle X^9 - 11 X^8 + 48 X^7 - 105 X^6 + 121 X^5 - 73 X^4 + 20 X^3 + 6 X^2 - 3 X + 1 \\
5/19 &&  \scriptstyle X^9 - 9 X^8 + 32 X^7 - 55 X^6 + 45 X^5 - 19 X^4 + 16 X^3 - 10 X^2 - 3 X - 1 \\
6/19 &&  \scriptstyle -X^9 + 7 X^8 - 20 X^7 + 25 X^6 - X^5 - 31 X^4 + 24 X^3 + 6 X^2 - 9 X - 1 \\
7/19 &&  \scriptstyle -X^9 + 5 X^8 - 12 X^7 + 15 X^6 - 9 X^5 - X^4 + 4 X^3 - 2 X^2 - X + 1 \\
8/19 &&  \scriptstyle X^9 - 3 X^8 + 8 X^7 - 13 X^6 + 17 X^5 - 17 X^4 + 12 X^3 - 6 X^2 + X + 1 \\
9/19 &&  \scriptstyle X^9 - X^8 + 8 X^7 - 7 X^6 + 21 X^5 - 15 X^4 + 20 X^3 - 10 X^2 + 5 X - 1 \\
1/20 &&  \scriptstyle (X^4 - 8 X^3 + 19 X^2 - 12 X + 1) (X^2 - 3 X + 1) (X^2 - 5 X + 5) (X - 2) \\
3/20 &&  \scriptstyle -(X^5 - 8 X^4 + 21 X^3 - 20 X^2 + 7 X - 2) (X^4 - 6 X^3 + 11 X^2 - 6 X - 1) \\
7/20 &&  \scriptstyle -(X^5 - 4 X^4 + 7 X^3 - 4 X^2 - X + 2) (X^4 - 2 X^3 + X^2 + 2 X - 1) \\
9/20 &&  \scriptstyle (X^4 + 3 X^2 + 1) (X^3 - X^2 + 3 X - 2) (X^2 - X + 1) \\
\end{align*}

\bibliographystyle{plain}
\bibliography{MF}

\end{document}